\pdfoutput=1
\documentclass[12pt]{article}
\usepackage[utf8]{inputenc}
\usepackage{color}
\usepackage[dvipsnames]{xcolor}
\usepackage{amsmath}
\usepackage{parskip}
\usepackage[showframe=false]{geometry}
\usepackage{changepage}
\usepackage{graphicx}
\usepackage{ dsfont }
\usepackage{diagbox}
\usepackage[refpage]{nomencl}
\usepackage{amsthm}
\usepackage[compact]{titlesec}         
\titlespacing{\section}{10pt}{10pt}{10pt} 
\AtBeginDocument{
  \setlength\abovedisplayskip{10pt}
  \setlength\belowdisplayskip{10pt}}
\usepackage{tikz}
\usetikzlibrary{positioning}
\usepackage{amssymb}
\usepackage{tikz}
\usepackage{cancel}
\usepackage{caption}
\usepackage{subcaption}
\usepackage{ textcomp }
\usepackage{ gensymb }
\usepackage{float}
\usepackage{rotating}
\usepackage{tabularx}
\usepackage{delarray}

\newcommand{\tpmod}[1]{{\@displayfalse\pmod{#1}}}

\usepackage{amsmath,xcolor}

\usepackage{mathtools}
\usepackage[shortlabels]{enumitem}

\usepackage[english]{babel}
\usepackage{soul}
\newtheorem{theorem}{Theorem}[section]

\theoremstyle{definition}

\usepackage{tikz}
\usepackage{siunitx}
\newcommand{\pj}{%
\begin{tikzpicture}%
\draw (0,-1.5pt) -- (8pt,-5pt);%
\draw (0,-1.5pt) -- (8pt,2pt);%
\draw (8pt,-5pt) -- (6pt,-6pt);%
\draw (8pt,-5pt) -- (6.5pt,-2.5pt);%
\end{tikzpicture}}
\usepackage{graphicx}
\title{Subtilings of Elnitsky Tilings for Finite Irreducible Coxeter Groups}
\author{Robert Nicolaides and Peter Rowley}
\date{2021}

\begin{document}

\maketitle
\begin{abstract}
Two new Elnitsky tilings for Coxeter groups of type $\mathrm{B}$ are displayed as certain subtilings. Additionally, a new tiling for the  non-crystallographic Coxeter group of type $\mathrm{H}_3$ is obtained, described as a $\mathrm{D}_6$-subtiling.
\end{abstract}
\section{Introduction}

In \cite{elnitsky} S. Elnitsky displays an elegant bijection between rhombic tilings of $2n$-gons and commutation classes of reduced words in the Coxeter group of type $\mathrm{A}_{n-1}$. Analogous results were also presented in \cite{elnitsky} for Coxeter groups of type $\mathrm{B}$ and $\mathrm{D}$. The main aim of this note is to make new tilings from old, with the aid of results from M\"uhlherr \cite{muhlherr}. Specifically, we uncover two new tilings for Coxeter groups of type $\mathrm{B}$, while also  giving an alternative (and easier) verification of one of Elnitsky tilings for type $\mathrm{B}$ described in Section 6 of \cite{elnitsky}. A new tiling for the non-crystallographic Coxeter group of type $\mathrm{H}_3$ is also produced. For recent work relating to tilings see  Escobar, Pechenik, Tenner,  Yong \cite{epty} and Tenner \cite{tenner1}, \cite{tenner2}. \\

Suppose $(W,S)$ is a Coxeter system. So $W$ is a group with $S$ a set of generators $S \subseteq W$ having a presentation
$$W = \langle S \;|\; (rs)^{m(r,s)} = 1\; \mbox{for} \; r,s \in S \rangle $$
where $m(s,s) = 1$, 
$m(r,s) = m(s,r) \geq 2$ for 
$s \ne r$ in $S$. Here we further assume that $W$ is a finite irreducible Coxeter group. In \cite{muhlherr} the concept of an admissable partition of $S$ is introduced -- we say more on this in Section \ref{background}. An admissable partition $\Sigma$ of $S$ gives rise to an embedding of a Coxeter system $(X,T)$ into $(W,S)$. This embedding, which we denote by $X\pj_{\Sigma} W$, has the important property that reduced words in $X$ correspond to reduced words in $W$. See Theorem \ref{muhlherr} for a precise statement of this.\\

Let $w \in W$. Then the length of $w$, denoted by $\ell(w)$, is the least $k$ such that $w = s_{i_1}s_{i_2} \dots s_{i_k}$ where $s_{i_j} \in S$. Any expression with $w = s_{i_1}s_{i_2} \dots s_{i_{\ell}}$ and $\ell = \ell(w)$ is called a reduced word for $w$. Now suppose that $(rs)^{m(r,s)} = 1$ with $r \ne s$ is one of the relations in the presentation for $W$. If a reduced word for $w$ has a consecutive subword $srsr \dots$, with a total of $m(r,s)$ $s$'s and $r$'s, then we can replace that subword by $rsrs \dots$, yielding another reduced word for $w$. Given a subset $J$ of the set of relations for $W$, we call two reduced words for $w$ $J$- equivalent if one can be obtained from the other by replacing a sequence of alternating subwords corresponding to relations in $J$. The resulting equivalence classes of the set of reduced words for $w$ will be denoted by $\mathcal{R}_{W,J}(w)$. We also denote the set of tilings associated with $w$ by $\mathcal{T}_W(w)$.\\

 \begin{theorem}\label{maintheorem}Suppose that $(W,S)$ is either of type $\mathrm{A}$ or type $\mathrm{D}$, and that $\Sigma = \{ \Sigma_i \; | \; i \in I \}$ is an admissable partition of $S$ with  $X\pj_{\Sigma} W$. Let $J$ be a subset of relations of $W$ to be described in Section \ref{background}. Put
 $$K = \{ (s_is_j)^{m(s_i,s_j)} \;|\;\mbox{for all} \; a\in \Sigma_i \; \mbox{and for all} \; b \in \Sigma_j, (s_as_b)^{m(s_a,s_b)} \in J  \} .$$
 For $x \in X$, let $\mathcal{T}_X(x)$ be the subset of $\mathcal{T}_W(x)$ corresponding to $\Sigma$-consistent reduced words of $x$. Then there exists as bijection between $\mathcal{R}_{X,K}(x)$ and $\mathcal{T}_X(x)$.
  \end{theorem}
  
This paper is arranged as follows. Section \ref{background} first expands on the work of M\"uhlherr already mentioned. Then we review in some detail the tilings for type $\mathrm{A}, \mathrm{B}$ and $\mathrm{D}$ which appear in \cite{elnitsky}. In Section \ref{subtilings} after proving Theorem \ref{maintheorem}, we then display a number of subtilings which yield tilings for Coxeter groups of type $\mathrm{B}$ and type $\mathrm{H}_3$.

\section{Background}\label{background}

We begin with the classification of admissable partitions for finite irreducible Coxeter groups. Let $(W,S)$ be a Coxeter system and $\Sigma = \{ \Sigma_i \; | \; i \in I \}$ a partition of $S$. If for each $i \in I$,  $W_{\Sigma_i} = \langle s \in S \; | \; s \in \Sigma_i \rangle$ is a finite parabolic subgroup of $W$, then we call 
$\Sigma$ spherical. For $\Sigma_i \in \Sigma$ we set $s_{\Sigma_i} = \Pi_{s \in \Sigma_i} s$, where some ordering of each $\Sigma_i$ is understood, though here we are only concerned with either the case when all the $s \in \Sigma_i$ pairwise commute or $s_{\Sigma_i}$ is the longest element of $W_{\Sigma_i}$. Set $T = \{s_{\Sigma_i} \; | \;\ i \in I \}$ and $X = \langle T \rangle$. Let $w \in W$. We say $w$ is $\Sigma$-consistent if it can be partitioned in the form $w = x_1x_2 \dots x_k$ where each $x_j$ is a reduced word for some $s_{\Sigma_i}$. We call $\Sigma$ admissable if for all $i \in I$ and for all $w \in W$ we have either $\Sigma_i \subseteq D_R(w)$ or $\Sigma_i \subseteq S \setminus D_R(w)$. Here $D_R(w) = \{ s \in S \; | \; \ell(w) = \ell(ws) - 1 \}$ denotes the right descent set of $w$ in the Coxeter system $(W,S)$.\\

For the finite irreducible Coxeter groups we consider their Coxeter/Dynkin diagrams are labelled as below, with $s_i$ or $t_i$ (depending on whether we are dealing with $W$ or $X$)  corresponding to the node labelled $i$.

\begin{figure}[H]
    \centering
\begin{tikzpicture}[scale=0.60,baseline={(0,0.0)}]
\node (s1) at (-2,0) {$\mathrm{A}_m$};
\draw (0,0) -- (4,0);

\filldraw [black] (5,-0) circle (1pt);
\filldraw [black] (6,-0) circle (1pt);
\filldraw [black] (7,-0) circle (1pt);
\draw(8,0) -- (10,0);
\node (s1) at (0,-0.6) {$1$};
\node (s1) at (2,-0.6) {$2$};
\node (s1) at (4,-0.6) {$3$};
\node (s1) at (10,-0.6) {$m$};

\filldraw [black] (0,-0) circle (5pt);
\filldraw [black] (2,-0) circle (5pt);
\filldraw [black] (4,-0) circle (5pt);
\filldraw [black] (8,-0) circle (5pt);
\filldraw [black] (10,-0) circle (5pt);

\node (s1) at (-2,0-4) {$\mathrm{B}_m$};
\draw (0,0-4) -- (4,0-4);
\node (s1) at (1,0.5-4) {$4$};

\filldraw [black] (5,-0-4) circle (1pt);
\filldraw [black] (6,-0-4) circle (1pt);
\filldraw [black] (7,-0-4) circle (1pt);
\draw(8,0-4) -- (10,0-4);
\node (s1) at (0,-0.6-4) {$1$};
\node (s1) at (2,-0.6-4) {$2$};
\node (s1) at (4,-0.6-4) {$3$};
\node (s1) at (10,-0.6-4) {$m$};

\filldraw [black] (0,-0-4) circle (5pt);
\filldraw [black] (2,-0-4) circle (5pt);
\filldraw [black] (4,-0-4) circle (5pt);
\filldraw [black] (8,-0-4) circle (5pt);
\filldraw [black] (10,-0-4) circle (5pt);

\node (s1) at (-2,0-4-4) {$\mathrm{D}_m$};
\draw (2,0-4-4) -- (4,0-4-4);
\draw (2,0-4-4) -- (0,1-4-4);
\draw (2,0-4-4) -- (0,-1-4-4);

\filldraw [black] (5,-0-4-4) circle (1pt);
\filldraw [black] (6,-0-4-4) circle (1pt);
\filldraw [black] (7,-0-4-4) circle (1pt);
\draw(8,0-4-4) -- (10,0-4-4);
\node (s1) at (0,1-4-4-0.6) {$1$};
\node (s1) at (0,-1-4-4-0.6) {$2$};
\node (s1) at (2,-0.6-4-4) {$3$};
\node (s1) at (4,-0.6-4-4) {$4$};
\node (s1) at (10,-0.6-4-4) {$m$};

\filldraw [black] (0,1-4-4) circle (5pt);
\filldraw [black] (0,-1-4-4) circle (5pt);
\filldraw [black] (2,-0-4-4) circle (5pt);
\filldraw [black] (4,-0-4-4) circle (5pt);
\filldraw [black] (8,-0-4-4) circle (5pt);
\filldraw [black] (10,-0-4-4) circle (5pt);

\node (s1) at (-2,0-4-4-4) {$\mathrm{E}_m$};
\node (s2) at (-2,0-4-4-4-1) {$(m=6,7,8)$};
\draw (0,0-4-4-4) -- (6,0-4-4-4);
\draw (4,0-4-4-4) -- (4,-2-4-4-4);

\filldraw [black] (9,-0-4-4-4) circle (1pt);
\filldraw [black] (8,-0-4-4-4) circle (1pt);
\filldraw [black] (7,-0-4-4-4) circle (1pt);
\node (s1) at (0,0.6-4-4-4) {$1$};
\node (s1) at (2,0.6-4-4-4) {$3$};
\node (s1) at (4,0.6-4-4-4) {$4$};
\node (s1) at (4,-2-0.6-4-4-4) {$2$};
\node (s1) at (6,0.6-4-4-4) {$5$};
\node (s1) at (10,0.6-4-4-4) {$m$};

\filldraw [black] (0,-0-4-4-4) circle (5pt);
\filldraw [black] (4,-2-0-4-4-4) circle (5pt);
\filldraw [black] (2,-0-4-4-4) circle (5pt);
\filldraw [black] (4,-0-4-4-4) circle (5pt);
\filldraw [black] (6,-0-4-4-4) circle (5pt);
\filldraw [black] (10,-0-4-4-4) circle (5pt);

\node (s1) at (-2,0-4-4-4-4) {$\mathrm{F}_4$};
\draw (0,0-4-4-4-4) -- (6,0-4-4-4-4);

\node (s1) at (0,-0.6-4-4-4-4) {$1$};
\node (s1) at (2,-0.6-4-4-4-4) {$2$};
\node (s1) at (4,-0.6-4-4-4-4) {$3$};
\node (s1) at (6,-0.6-4-4-4-4) {$4$};

\node (s1) at (3,0.5-4-4-4-4) {$4$};

\filldraw [black] (0,-0-4-4-4-4) circle (5pt);
\filldraw [black] (2,-0-4-4-4-4) circle (5pt);
\filldraw [black] (4,-0-4-4-4-4) circle (5pt);
\filldraw [black] (6,-0-4-4-4-4) circle (5pt);

\node (s1) at (-2,0-4-4-4-4-4) {$\mathrm{H}_3$};
\draw (0,0-4-4-4-4-4) -- (4,0-4-4-4-4-4);

\node (s1) at (0,-0.6-4-4-4-4-4) {$1$};
\node (s1) at (2,-0.6-4-4-4-4-4) {$2$};
\node (s1) at (4,-0.6-4-4-4-4-4) {$3$};

\node (s1) at (1,0.5-4-4-4-4-4) {$5$};

\filldraw [black] (0,-0-4-4-4-4-4) circle (5pt);
\filldraw [black] (2,-0-4-4-4-4-4) circle (5pt);
\filldraw [black] (4,-0-4-4-4-4-4) circle (5pt);

\node (s1) at (-2+8,0-4-4-4-4-4) {$\mathrm{H}_4$};
\draw (0+8,0-4-4-4-4-4) -- (6+8,0-4-4-4-4-4);

\node (s1) at (0+8,-0.6-4-4-4-4-4) {$1$};
\node (s1) at (2+8,-0.6-4-4-4-4-4) {$2$};
\node (s1) at (4+8,-0.6-4-4-4-4-4) {$3$};
\node (s1) at (6+8,-0.6-4-4-4-4-4) {$4$};

\node (s1) at (1+8,0.5-4-4-4-4-4) {$5$};

\filldraw [black] (0+8,-0-4-4-4-4-4) circle (5pt);
\filldraw [black] (2+8,-0-4-4-4-4-4) circle (5pt);
\filldraw [black] (4+8,-0-4-4-4-4-4) circle (5pt);
\filldraw [black] (6+8,-0-4-4-4-4-4) circle (5pt);

\end{tikzpicture}
    \caption{The finite irreducible Coxeter groups excluding the dihedral groups.}
    \label{Figure FICG}
\end{figure}

 We now present the result on admissible embeddings we need for our subtilings.
 
 \begin{theorem} (M\"uhlherr \cite{muhlherr})\label{muhlherr} Suppose that $W$ is a finite irreducible Coxeter group.
 \begin{enumerate}
 \item[(i)] The admissable partitions $\Sigma$ giving  $X \pj_{\Sigma} W$ are listed in Table \ref{Table Admissible Partitions} (where the cases of either $X$ or $W$ being a dihedral group has been omitted).
 \item [(ii)]If $X \pj_{\Sigma} W$, then  $x$ is 
 $\Sigma$-consistent for all $x \in X$.
 \end{enumerate}
\end{theorem}

\begin{table}[H]
    \centering
    \begin{tabular}{c c c}
    \vspace{+0.2cm} 
        Type of $W$ & Type of $X$ & $\Sigma$\\\vspace{+0.2cm} 
        $\mathrm{A}_{2n-1} \, (n\ge 2)$ & $\mathrm{B}_n$ & $\{\{i,2n-i\},\{n\} \; | \; i=1,\ldots,n-1\}$\\\vspace{+0.2cm} 
        $\mathrm{A}_{2n}  \, (n\ge 2)$ & $\mathrm{B}_n$ & $\{\{i,2n-i+1\} \;|  \;i=1,\ldots,n\}$\\\vspace{+0.2cm} 
        $\mathrm{D}_{n+1}  \, (n\ge 2)$ & $\mathrm{B}_n$ & $\{\{1,2\},\{i\} \; | \; i=3,\ldots,n+1\}$\\\vspace{+0.2cm} 
        $\mathrm{E}_6 $ & $\mathrm{F}_4$ & $\{\{1,6\},\{3,5\},\{2\},\{4\}\}$\\\vspace{+0.2cm} 
        $\mathrm{D}_6 $ & $\mathrm{H}_3$ & $\{\{1,4\},\{2,6\},\{3,5\}\}$\\\vspace{+0.2cm} 
        $\mathrm{E}_8 $ & $\mathrm{H}_4$ & $\{\{1,8\},\{2,5\},\{3,7\},\{4,6\}\}$\\
    \end{tabular}
    \caption{Admissible partitions for the finite irreducible Coxeter groups} 
    \label{Table Admissible Partitions}
\end{table}
 We next outline Elnitsky's original construction for tilings. Consider $(W,S) = (\mathrm{Sym}(n),\{(i,i+1) \; | \; i \in \{1, \dots n-1\} \})$, the type $\mathrm{A}$ Coxeter system of rank $n-1$. Let $w \in W$ be some permutation, then we construct $Y(w)$, a (possibly degenerate) $2n$-gon. First we describe a \emph{regular} construction of $Y(w)$ before mentioning the degrees of freedom we can employ to preserve the properties we care about.

\begin{enumerate}[$(i)$]
    \item Let $Y(w)$ be a $2n$-gon with unit side lengths and whose upper-most vertex we call $M$.
    \item Construct and label the first $n$ edges anti-clockwise from $M$ with labels from $\{1, \dots, n\}$ consecutively, ensuring that they form half of a regular $2n$-gon, so as the angles between edges are the same. 
    \item  Construct and label the first $n$ edges clockwise from $M$ such that the $i^{th}$ edge from $M$ is parallel to $w^{-1}(i)$. 
\end{enumerate}

The choice of angles and side lengths for the left hand side do not actually matter so long as we have the left hand side's vertices forming a convex set and the edges with the same label being translations of one another. We will keep to the regular case in our examples where possible and call the right-most path in $Y(w)$ the border of $w$ and denote it $B(w)$. This naturally gives us a bijection between $W = \mathrm{S}(n)$ and $\{B(w) \; | \; w \in W\}$. For example, the left-most path for all $X(w)$ is always $B(id)$.\\

We now consider \emph{tilings} of $Y(w)$ by rhombi. Implicitly, by rhombic tiling we mean a covering of $Y(w)$ by regions of rhombi that intersect only on their boundaries. Denote the set of all rhombic tilings of $Y(w)$ by $\mathcal{T}_W(w)$. 

Let $J = \{(s_is_j)^2 \;| \; |i-j|\ge 2 \}$, a subset of the relations for $W$. So $\mathcal{R}_{W,J}(w)$, is the set of reduced words of $w$ over $S$ up to commuting generators. Elnitsky is able to prove quite directly the following elegant fact.

\begin{theorem}\label{Elnitsky A}
For all $w \in W = \mathrm{S}(n)$, there exists a bijection between $\mathcal{T}_W(w)$ and $\mathcal{R}_{W,J}(w)$.
\end{theorem}

In \cite{elnitsky}, we also learn of similar constructions for types $\mathrm{B}$ and $\mathrm{D}$. Those of type $\mathrm{B}$ are exactly the tilings of $\mathrm{A}$ that are \emph{horizontally symmetric} - they can be flipped about the horizontal line through the middle vertex. It is essentially the same as that of type $\mathrm{A}$ and we will re-examine it in 3.1. We now remark on the construction for type $\mathrm{D}$ Coxeter groups. 
Let $(W,S)$ be $(\mathrm{D}_n,\{s_1,s_2,\ldots,s_n\})$, a standard embedding of $\mathrm{D}_n$ into $\mathrm{Sym}(2n)$ with $s_1=(1,-2)(2,-1), s_2 =(2,3)(-2,-3), s_3=(1,2)(-1,-2)$ and $s_i = (i-1,i)(-(i-1),-i)$ for each $i \in \{4, \dots, n \}$. In this setting, Elnitsky describes a similar construction for tilings of a $4n$-gon such that we again have a correspondence of tilings and reduced words of the Coxeter group up to some commuting generators. 

Here is the construction of the polygon $Y(w)$ for all $w \in W$.
\begin{enumerate}[$(i)$]
    \item Let $U$ be the upper most vertex of our $4n$-gon, and $L$ the lower most vertex and $M$ the vertex that is an equal distance from both.
    \item Let the first $2n$ edges anticlockwise from $U$ be those of the regular $4n$-gon with unit length edges, labelling them with $\{n,n-1,\ldots,1,-1,\ldots,-n\}$ respectively.
    \item  For each $i \in \{-n,-n-1,\ldots,-1,1,\ldots,n\}$, construct and label the $2n$ edges anti-clockwise from $L$ such that the $i^{th}$ edge from $M$ in this direction has the same length as $w^{-1}(i)$ and parallel be to it. 
\end{enumerate}
There is one extra non-trivial condition that is necessary for the desired correspondence, namely that the absolute gradient of each of our edges from the horizontal is always $> \pi/3$. It is a non-trivial observation by Elnitsky that this removes certain avoidable intersections of tiles. \\

Let us reuse the language of $B(w)$ to denote the rightmost path of $Y(w)$. This time we consider a new set of tiling rules which are more complex than that of type $\mathrm{A}$. In particular we now have a set of \emph{megatiles} at our disposal. The megatiles are a subset of octagons with unit edge lengths whose construction we now discuss. Its upper-most vertex, $U$, and lower-most vertex, $L$, must lie on a vertical line. Its first four edges anti-clockwise from $U$ must be symmetric through the horizontal line passing through the middle vertex. Call these edges $E$. Then to make the remaining edges perform the following on $E$.
\begin{enumerate}
    \item[(i)] Transpose the first and second pair of edges and the third and fourth pair of edges in $E$.
    \item[(ii)] Reflect $E$ through the vertical line passing through $U$ and $L$. 
\end{enumerate}

We call the set of all such tilings for $Y(w)$, $\mathcal{T}_W(w)$. Then Elnitsky proves the following by a direct argument where the relation set $J$ consists of $(s_is_j)^2$ for $i,j \in \{1, \dots, n \}$ with $|i-j|\ge 2$ but excluding $(s_1s_3)^2$. 
\begin{theorem}\label{Elnitsky D}
For all $w \in W$ with $W$ of type $\mathrm{D}_n$, there is a bijection between $\mathcal{T}_W(w)$ and $\mathcal{R}_{W,J}(w).$
\end{theorem}

\section{Subtilings for type $\mathrm{B}$ and type $\mathrm{H}_3$}\label{subtilings}

First we prove Theorem \ref{maintheorem}.

$\mathbf{Proof \; of \;Theorem}$ \ref{maintheorem}
By Theorems \ref{Elnitsky A} and \ref{Elnitsky D} there exists a bijection between $\mathcal{R}_{W,J}(w)$ and $\mathcal{T}_W(w)$. For an admissible partition  $\Sigma$ of $S$, Theorem \ref{muhlherr}(ii) says all the reduced $\Sigma$-consistent words are reduced in $W$. These are in correspondence with the reduced words of $X$. Therefore, we have the bijection between the restricted set of words and tiles. To see that the associated relation set must be $K$ we observe that for $\{\Sigma_i,\Sigma_j \}\subseteq \Sigma$ to be in the relation set, we must have that for all elements of $\Sigma_i$ and all elements of $\Sigma_j$ they must be relations in $J$ . 



Combining Theorems \ref{maintheorem} and \ref{muhlherr} yields subtilings which we now discuss. For  3.1 to 3.4 we assume $X \pj_{\Sigma} W$ with Coxeter systems $(W,S)$ and $(X,T)$.  Let $x_0$ be the longest element of $X$.

\subsection{$X$ of type $\mathrm{B}_n$, $W$ of type $\mathrm{A}_{2n-1}$}\label{3.1}
So for our subtiling we are implicitly using the embedding, $t_1 \rightarrow s_n=(n,n+1)$ while $t_i \rightarrow s_{n-i}s_{n+i}=(n-i,n-i+1)(n+i,n+i+1)$ for $i \in \{1, \dots, n-1 \}$. In this case, the relation set $K$ is $\{(t_it_j)^2 \; | \; |i-j|\ge 2\}$.
This gives us the tilings of type $\mathrm{A}$ that are symmetric about their horizontal axis. These are exactly Elnitsky's tilings of type $\mathrm{B}$ in Section 6 of \cite{elnitsky}. Here we have an alternative proof of this without the gritty geometric details appearing in \cite{elnitsky}.

We demonstrate some of the examples for $\mathcal{T}_X(x_0)$ 
when $n=3$.

\begin{figure}[H]   \centering   \includegraphics[width = 15cm]{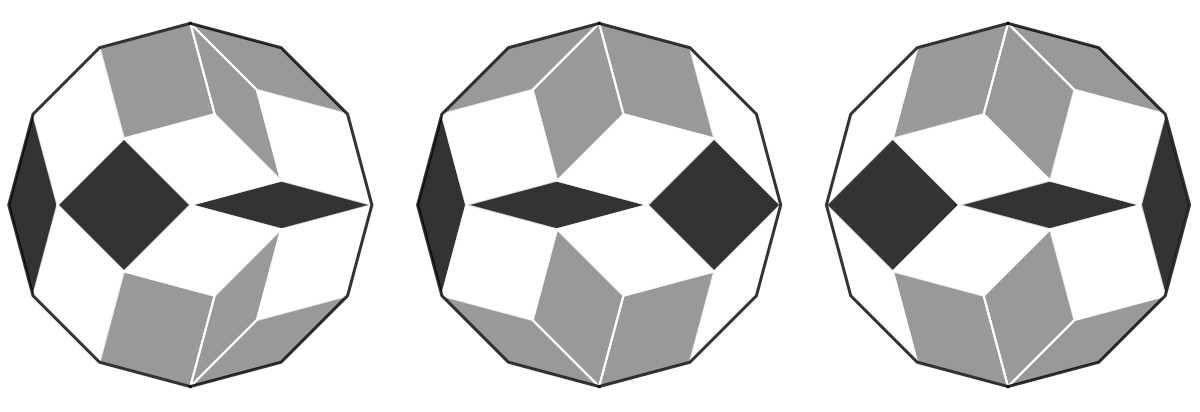}    \caption{Elnitsky tiling of type $\mathrm{B_3}$ viewed as a subtiling of the type $\mathrm{A_5}$ tiling.}   \label{Etiling B1}\end{figure}

\subsection{$X$ of type $\mathrm{B}_n$, $W$ of type $\mathrm{A}_{2n}$}
\begin{figure}[H]
    \centering
    \includegraphics[width = 15cm]{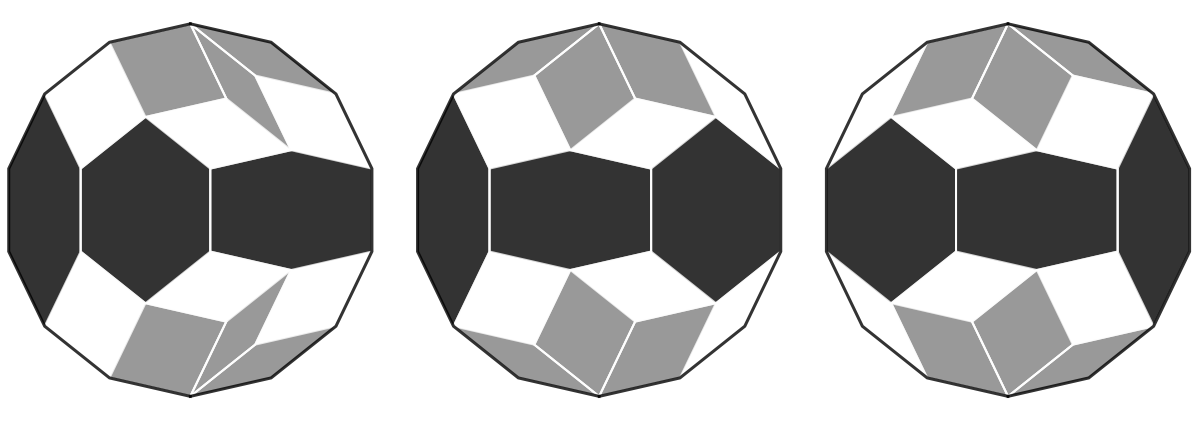}
    \caption{Elnitsky tiling of type $\mathrm{B}_3$ viewed as a subtiling of the type $\mathrm{A}_6$ tiling.}
    \label{Etiling B2}
\end{figure}

The other admissible partition of type $\mathrm{A}$ groups that induce a Coxeter group of type $\mathrm{B}$ is the following.
Here $t_1$ is sent to the longest element in the parabolic subgroup of $\{s_n,s_{n+1}\}$, which is $t_1 \rightarrow s_ns_{n+1}s_n = s_{n+1}s_ns_{n+1}=(n,n+2)$, while the others are sent to $t_i \rightarrow s_{n-i}s_{n+1+i} = (n-i,n-i+1)(n+1+i,n+2+i)$.  Since $t_1$ is not the product of disjoint transpositions, the corresponding tile is necessarily formed by placing the sequence of tiles corresponding to either $s_1s_2s_1$or $s_2s_1s_2$ in Elnitstky’s type $\mathrm{A}$ tiling. We identifying the placement of these equivalent sequences as one so-called megatile which itself is a hexagon. Again, we observe that $K =  \{(t_it_j)^2 \; |\;  |i-j|\ge 2\}$ and we again consider what $\mathcal{T}_X(x_0)$ looks like for the case $n=3$ -- see Figure \ref{Etiling B2}.

We observe that the existence of this particular tiling given that of Elnitsky’s type $\mathrm{B}$ tiling is, in hindsight, very intuitive; as it has horizontal symmetry, if we were to insert the constant vertical edge in place of the middle vertex we will preserve the tiling and words. Similar observations can be found in \cite{elnitsky, hamanaka}, when studying strips. 

\subsection{$X$ of type $\mathrm{B}_n$, $W$ of type $\mathrm{D}_{n + 1}$}
We consider the final admissible partition that induces $\mathrm{B}_n$. This time it is a partition of $\mathrm{D}_{n+1}$.

In this case we have $t_1 \rightarrow s_1s_2 = (n-1,n+1)(n,n+2)(1,2)(2n-1,2n)$ along with $t_i \rightarrow s_{i + 1} = (i,i+1)(2n-i,2n-i-1)$ for $2\le i \le n-2$.  We look at the reduced words of longest element of $B_3$ this time. Note that this time the relation set $K$ is $\{ (t_it_j)^2 \;| \; |i-j|\ge 2\}\setminus\{ (t_1t_3)^2\}$. In particular, this means that when $n=4$ we have $K$ is empty and consequently we now get six tiles corresponding to $x_0$.

\begin{figure}[H]
    \centering
    \includegraphics[width = 15cm]{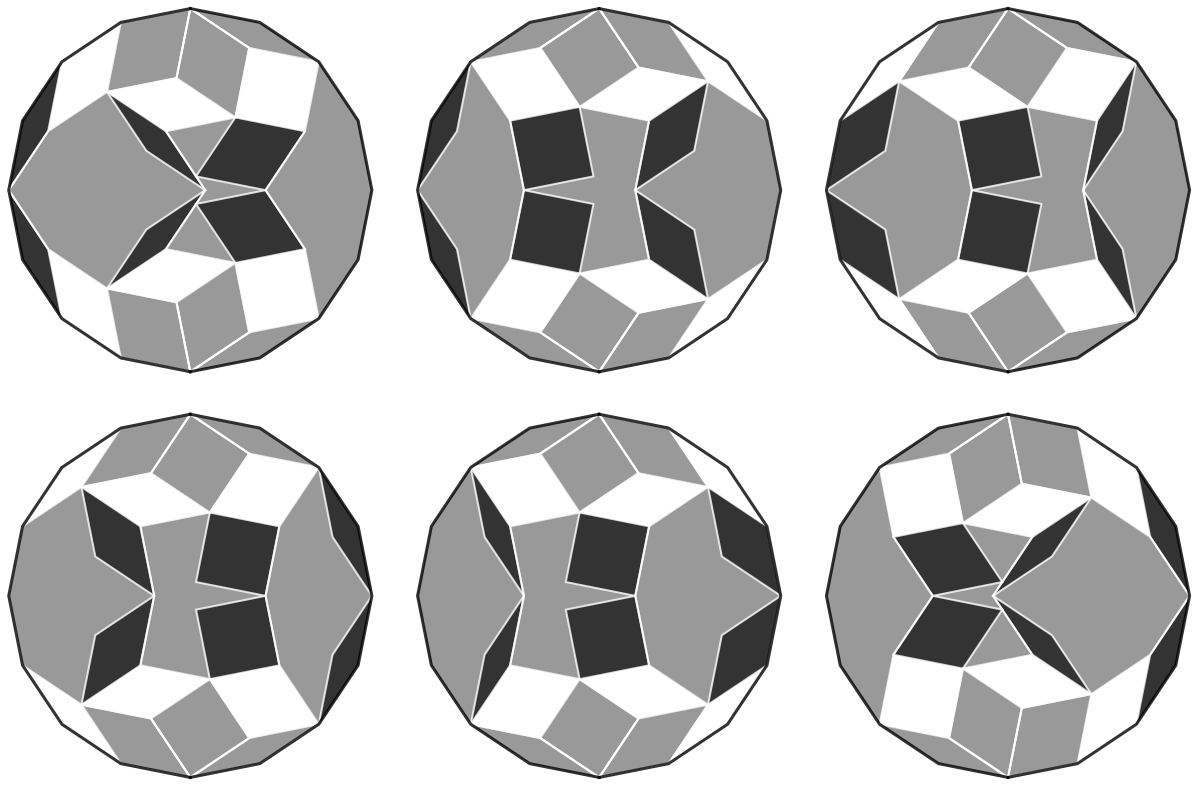}
    \caption{Elnitsky tiling of type $\mathrm{B}_4$ viewed as a subtiling of the type $\mathrm{D}_5$ tiling.}
    \label{Etiling B in D}
\end{figure}

We have chosen a regular polygon to be the border of this element despite this producing the fixable intersections mentioned in Section 2.

\subsection{$X$ of type $\mathrm{H}_3$, $W$ of type $\mathrm{D}_6$}
Finally, we consider the tiling for $\mathrm{H}_3$ as a subtiling for $\mathrm{D}_6$ induced from the following admissible partition.

\begin{figure}[H]
\centering

\begin{tikzpicture}[baseline={(0,0.0)}]
\node (d1) at (1,-2.5) {$s_2$};
\node (d4) at (1,2.5) {$s_{1}$};
\node (d2) at (4,-1) {$s_3$};
\node (d2) at (6,-1) {$s_4$};
\node (d2) at (8,-1) {$s_5$};
\node (d3) at (10,-1) {$s_{6}$};

\draw (1,2) -- (4,0);
\draw (1,-2) -- (4,0);
\draw (4,0) -- (6,0);
\draw (8,0) -- (10,0);
\draw  (6,0) -- (8,0);

\filldraw [black] (1,2) circle (5pt);
\filldraw [gray] (1,-2) circle (5pt);
\filldraw [white, draw=black] (4,-0) circle (5pt);
\filldraw [gray] (6,0) circle (5pt);
\filldraw [white, draw=black] (8,-0) circle (5pt);
\filldraw [black] (10,-0) circle (5pt);

\draw[->] (6.5+0.5,-2) -- (6.5+0.5,-3);

\node (d3) at (4+0.5+0.5,-5) {$t_{1}$};
\node (d4) at (6+0.5+0.5,-5) {$t_{2}$};
\node (d4) at (8+0.5+0.5,-5) {$t_{3}$};
\node (d4) at (5+0.5+0.5,-3.5) {$5$};

\draw (4+0.5+0.5,-4) -- (6+0.5+0.5,-4);
\draw (6+0.5+0.5,-4) -- (6.5+0.5+0.5,-4);
\draw (7.5+1,-4) -- (8+1,-4);
\draw (6.5+1,-4) -- (7.5+1,-4);

\filldraw [gray] (4+1,-4) circle (5pt);
\filldraw [white, draw=black] (6+1,-4) circle (5pt);
\filldraw [black] (8+1,-4) circle (5pt);

\end{tikzpicture}
\caption{The admissible partition of $\mathrm{D}_{6}$ into $\mathrm{H}_3$.}
\label{admiss D6 to H3}
\end{figure}
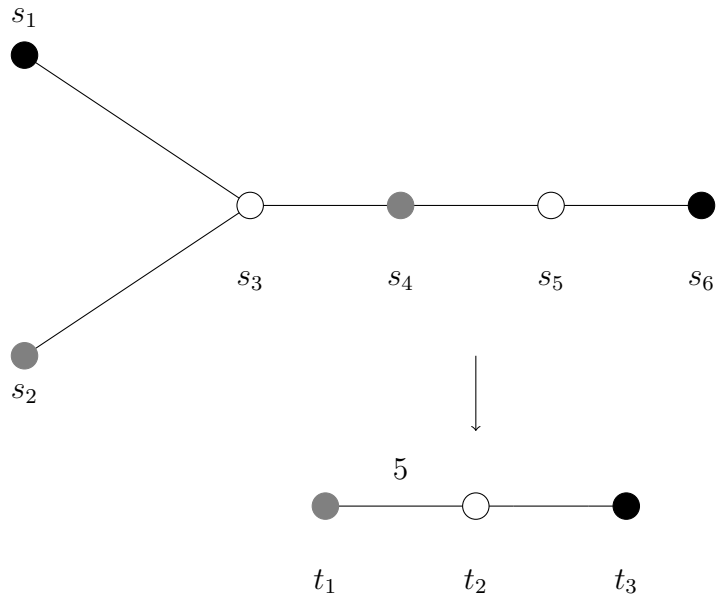

\begin{figure}[H]
    \centering
    \includegraphics[width = 15cm]{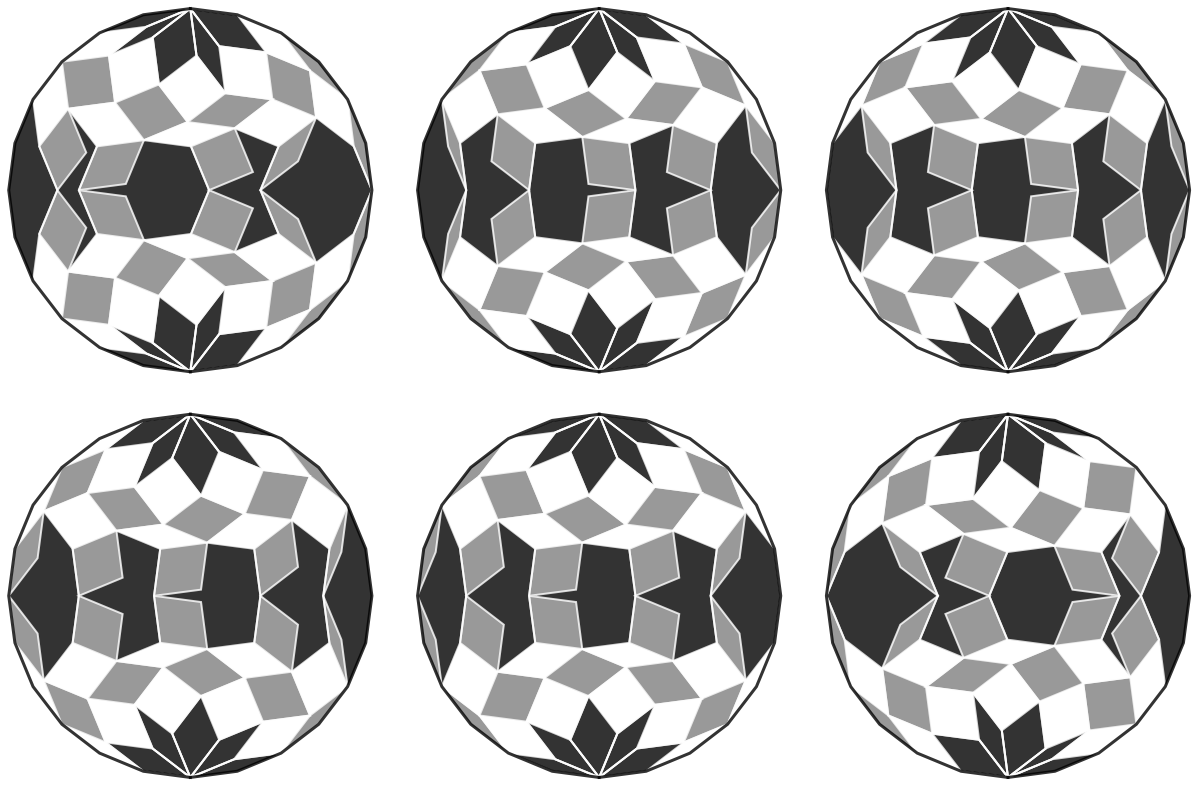}
    \caption{Elnitsky tiling of type $\mathrm{H}_3$ viewed as a subtiling of the type $\mathrm{D}_6$ tiling.}
    \label{Etiling H3}
\end{figure}

For this tiling we have an empty relation set, which, as luck would have it, gives us genuine bijections between reduced words $w$ of $\mathrm{H}_3$ and subtilings in $\mathcal{T}_X(w)$. There are 286 reduced words for the longest element of $\mathrm{H}_3$, we highlight a selected sample of six corresponding tilings in Figure \ref{Etiling H3}.

\end{document}